\documentclass[12pt,a4]{article}

\topmargin= - 10 true mm
\usepackage{amsmath,graphicx}
\usepackage{amsthm}
\usepackage{amssymb}
\usepackage{bbm}
\usepackage{setspace}

\newcommand{\exn}{{\bf E}}%{\mbox{\boldmath$E$}}
\newcommand{\pr}{{\bf P}}%{\mbox{\boldmath$P$}}
\newcommand{\N}{\mathbb{N}}
\newcommand{\R}{\mathbb{R}}

\newcommand{\deq}{\stackrel{d}{=}}
\newcommand{\ep}{\varepsilon}

\newtheorem*{theorem*}{Theorem}
\newcommand{\mf}[1]{\mathbf{#1}}
\newcommand{\lag}{\langle}  % left vector product bracket
\newcommand{\rag}{\rangle}  % right vector product bracket
\newcommand{\bla}{\mbox{\boldmath$\lambda$}}     % bold la
\newcommand{\bx}{\mbox{\boldmath$x$}}        %bold x
\newcommand{\bbxi}{\mbox{\boldmath$\xi$}}         % bold \xi
\newcommand{\bal}{\mbox{\boldmath$\alpha$}}   % bold \alpha
\newcommand{\sbx}{\mbox{\scriptsize\boldmath$x$}}     % small bold x
\newcommand{\sbla}{\mbox{\scriptsize\boldmath$\lambda$}}   % small bold \lambda
\newcommand{\sbal}{\mbox{\scriptsize\boldmath$\alpha$}}   % small bold \alpha
\newcommand{\bv}[1]{\mbox{\boldmath$#1$}}   % bold vector
\newcommand{\sbv}[1]{\mbox{\scriptsize\boldmath$#1$}}  % small bold vector
\DeclareMathOperator*{\argminA}{arg\,min}

\begin{document}

\title{The exact asymptotics for hitting probability of  a remote orthant by a multivariate L\'evy process: the Cram\'er case}

\author{K.~Borovkov$^1$ and Z.~Palmowski$^2$}

\date{}

\footnotetext[1]{School of Mathematics and Statistics, The University of Melbourne, Parkville 3010, Australia. E-mail: borovkov@unimelb.edu.au.}

\footnotetext[2]{Department of Applied Mathematics, Wroc{\l}aw University of Science and Technology,  27 Wybrze\.{z}e Wyspia\'nskiego st.,
50-370 Wroc{\l}aw, Poland.  E-mail: zbigniew.palmowski@gmail.com.}

%\date{}
\maketitle

\begin{abstract}
For a multivariate L\'evy process satisfying the Cram\'er moment condition and having a drift vector with at least one negative component, we  derive the exact asymptotics of the probability of ever hitting the positive orthant that is being translated to infinity along a fixed vector with positive components. This problem is motivated by the multivariate ruin problem introduced in F.~Avram et~al.\ (2008) in the two-dimensional case. Our solution relies on the analysis from Y.~Pan and K.~Borovkov~(2017) for multivariate random walks  and an appropriate time discretization.


\end{abstract}

\section{Introduction}
In this note we consider the following large deviation problem for continuous time processes with independent increments  that was motivated by the multivariate simultaneous ruin problem introduced in~\cite{pal1}.
Let $\{\bv{X}(t)\}_{t\ge 0}$ be a $d$-dimensional ($d\ge 2$) right-continuous  L\'evy process with $\bv{X}(0)=0$.
%and, for a Borel set $V\subset \R^d$, denote by
%\[
%  \tau (V) :=\inf\{t\ge 0: \bv{X}(t)\in V\}
%\]
%the first hitting time of the set $V$ by $\bv{X}$.
One is interested in finding  the precise asymptotics for the hitting probability
%$\pr (\tau  (sG)<\infty)$
of the   orthant $sG$ as $s\to\infty$, where
$$
G:= \bv{g}+ Q^+
$$
for some fixed
\[
\bv{g}\in   Q^+,
\quad
  Q^+ := \{
\bv{x}=(x_1, \ldots, x_d) \in \mathbb{R}^d: x_j > 0,\ 1 \le j \le d \}.
\]
Clearly, $sG= s\bv{g}+ Q^+,$ which is just the positive orthant translated by~$s\bv{g}.$

We solve this problem  under appropriate  Cram\'er moment assumptions and further  conditions on the process~$\bv{X}$ and vertex~$\bv{g} $ that, roughly speaking, ensure that the ``most likely place" for $\bv{X}$ to hit $sG$ when $s$ is large is in vicinity of the ``corner point"~$s\bv{g}$. More specifically, we show that the precise asymptotics of the  hitting probability of $sG$   are given by the following expression: letting
\[
  \tau (V) :=\inf\{t\ge 0: \bv{X}(t)\in V \}
\]
be the first hitting time of the set $V\subset \R^d$ by the process~$\bv{X},$ one has
\begin{equation}
\label{first}
\mf{P}\big( \tau(sG) < \infty \big) =
A_0 s^{-(d-1)/2}e^{-sD(G)}(1+o(1))
\quad \mbox{as}\quad  s\to\infty,
\end{equation}
where the ``adjustment coefficient" $D(G)$ is the value of the second rate function (see~\eqref{Dv} below) for the distribution of $\bv{X}(1)$  on the set~$G$ and the constant $A_0 \in (0, \infty)$ can be computed explicitly.

The asymptotics \eqref{first} extend a number of  known results.
The main body of literature on the topic of precise asymptotics for boundary crossing large deviation probabilities in the multivariate case concerns the random walk theory, see \cite{paper23,Borbook,mainpaper} and references therein for an overview of  the  relevant results. The crude logarithmic asymptotics in the multivariate case was also derived independently in~\cite{Coll}.

The entrance probability to a remote set for L\'evy processes was analyzed later, usually under some specific assumptions on the structure of these processes. For example,  paper~\cite{pal1} dealt with   the two-dimensional reserve process of the form
\begin{equation}
\label{CL}
\bv{X}(t)=(X_1(t), X_2(t))=(c_1,c_2)\sum_{i=1}^{N(t)}C_i-(p_1,p_2)t,\quad t\ge 0,
\end{equation}
where $c_i, p_i>0,$ $i=1,2,$ are constants, $\{C_i\}_{i\ge 1}$ is a sequence of  i.i.d.\  claim sizes, and $N(t)$ is an independent of the claim sizes Poisson process. That model admits the following interpretation: the  process $\bv{X}$ describes the dynamics of the reserves $X_i(t),$  $i=1,2,$  of two insurance companies that divide between them both claims and premia in some pre-specified proportions. In that  case, $\mf{P}\big( \tau(sG) < \infty \big)$ corresponds to the simultaneous ruin probability of the two companies.  The main result of the present paper generalizes the assertion of Theorem~5 of~\cite{pal1} to the case of general L\'evy processes. One may also wish to mention here the relevant  papers~\cite{BD,PalPist}.

\section{The main result}
To state the main result,  we will need  some notations. For brevity, denote by $\bv{\xi}$ a random vector such that
\begin{equation}\label{defxi}\bv{\xi}\deq\bv{X}(1).
\end{equation}
Our first condition on~$\bv{X}$ is stated as follows.

\smallskip
 [$\textbf{C}_1$]~{\em The distribution of $\bbxi$ is non-lattice and there is no hyperplane $H = \{ \bx :  \lag \bv{a}, \bx \rag = c \} \subset \R^d$ such that $\pr (\xi\in H) = 1.$}
\smallskip

That condition can clearly  be re-stated in terms of the covariance matrix of the Brownian component and   spectral measure of $\bv{X}$, although such re-statement will not make it more compact nor transparent.

Next denote by
\begin{equation}\label{cumulant}
K(\bla):=\ln \exn  e^{ \langle \sbla , \sbv{\xi}\rangle}, \quad \bv{\lambda} \in \R^d,
\end{equation}
the cumulant function of~$\bv{\xi}$ and let
\begin{equation*}
\Theta_{\psi} := \{ \bv{\lambda} \in \R^d:  K(\bv{\lambda}) < \infty \}
\end{equation*}
be the set on which the moment generating function of $\bv{\xi}$ is finite. We will need the following Cram\'{e}r moment condition on $\bv{X}$:

\smallskip
[$\textbf{C}_2$]~\textit{$\Theta_{\psi}$ contains a non-empty open set}.
\smallskip

The  first rate function $\Lambda(\bv{\alpha})$ for the random vector $\bv{\xi}$ is defined as the Legendre transform of the cumulant function~$K:$
\begin{equation}
\label{rf1}
   \Lambda(\bv{\alpha}) :=   \sup_{\sbv{\lambda} \in \Theta_{\psi}} (\langle \bv{\alpha},\bv{\lambda} \rangle - K(\bv{\lambda})), \quad  \bv{\alpha} \in \R^d.
\end{equation}
The probabilistic interpretation of the first rate function is given by the following relation (see e.g.~\cite{mainpaper}): for any $\bal \in \R^d$,
\begin{equation}
\label{rf2}
\Lambda(\bal) = -\lim_{\ep \to  0}\lim_{n \to \infty}
  \frac{1}{n}\ln \mf{P}\biggl(\frac{\bv{X}(n)}{n} \in U_{\ep}(\bal) \biggr),
\end{equation}
where $U_{\varepsilon}(\bal)$ is the $\varepsilon$-neighborhood of $\bal$.
Accordingly, for a set $B \subset \R^d$,  any point $\bal \in B$ such that
\begin{equation}
\label{MPPdef}
\Lambda(\bal) =\Lambda (B):= \inf_{\sbv{v} \in B}\Lambda(\bv{v})
\end{equation}
is called a most probable point   (MPP)  of the set $B$ (cf.\ relation~(11) in~\cite{PaBo}).
If such a point $\bal$ is unique for a given set~$B$, we denote it by
\begin{equation}
\label{MPP}
\bal[B] := \argminA_{\sbv{v} \in B}\Lambda(\bv{v}).
\end{equation}

Now recall the definition of the second rate function $D$ that was introduced and studied  in~\cite{paper23}: letting $ D_u(\bv{v}):=u\Lambda( \bv{v}/u)     $ for $\bv{v} \in \R^d,$ one sets
\begin{equation}
 \label{Dv}
D(\bv{v}) := \inf_{u > 0} D_u(\bv{v}), \quad  \bv{v}\in \R^d, \qquad
D(B):= \inf_{\sbv{v}\in B} D(\bv{v}), \quad B\subset \R^d
\end{equation}
(see also~\cite{PaBo}). Further, we put
\begin{equation}
\label{mpt}
r_B := \argminA_{r > 0} D_{1/r}(B).
\end{equation}
Recall the probabilistic meaning of the function~$D$ and the value $r_B$.
While the first rate function $\Lambda$ specifies the main term in  the asymptotics of the probabilities for the random walk values $\bv{X}(n)$ to be inside ``remote sets" (roughly speaking, $\Lambda (B)$ equals the RHS of~\eqref{rf2} with the neighbourhood of $\bal$ in it replaced with~$B$), the second rate function $D$ does that for the probabilities of {\em ever hitting\/} ``remote sets" by the whole random walk trajectory $\{\bv{X}(n)\}_{n\ge 0},$ the meaning of $r_B$ being that $1/r_B$ gives (after appropriate scaling) the ``most probable time" for the walk to hit the respective remote set. For more detail, we refer the interested reader to~\cite{paper23, PaBo}.

Define the Cram\'{e}r range $\Omega_{\Lambda}$ for $\bv{\xi}$ as follows:
\begin{equation*}
\Omega_{\Lambda}
 := \big\{ \bal ={\rm grad}\, K (\bla) : \bla \in {\rm int}(\Theta_{\psi})  \big\},
\end{equation*}
where the cumulant function $K(\bla)$ of $\bv{\xi}$ was defined in~\eqref{cumulant} and  $ {\rm int}(B)$ stands for the interior of the set~$B$. In words, the set $\Omega_{\Lambda}$ consists of all the vectors that can be obtained as the expectations of the Cram\'{e}r transforms of the law of $\bv{\xi}$, i.e.\  the distributions  of the form $  e^{ \langle \sbla, \sbx\rangle- K(\sbla )}\pr (\bv{\xi}\in d\bv{x})$, for parameter values $\bla \in {\rm int}(\Theta_{\psi})$.

For $\bal \in \R^d$, denote by $\bla(\bal)$ the vector $\bla$ at which the upper bound in \eqref{rf1} is attained (when such a vector exists, in which case it is always unique):
\begin{equation*}
 \Lambda(\bv{\alpha})
  =\langle \bv{\alpha},\bv{\lambda}(\bal)\rangle - K(\bv{\lambda}(\bal)).
\end{equation*}

For   $r>0,$ assuming that  $\bal[rG] \in \Omega_{\Lambda},$ introduce the vector
\begin{equation}
\label{normaln}
 \bv{N}(r) :=  {\rm grad}\, \Lambda(\bal)\big|_{\sbal = \sbal[rG]}= \bla(\bal[rG]),
\end{equation}
which is a normal   to the  level surface of $\Lambda$  at the point $\bal[rG]$
(see e.g.\ (22) in~\cite{PaBo}).

The last condition that we will need to state  our main result depends on the parameter $r > 0$ and is formulated  as follows:

\smallskip
[\textbf{C}$_3(r)$]~\textit{One has }
\begin{equation*}
\Lambda(rG) = \Lambda(r\bv{g}), \quad r\bv{g} \in \Omega_{\Lambda}, \quad  \bv{N}(r) \in Q^+, \quad \lag \mf{E} \bbxi, \bv{N}(r) \rag < 0.
\end{equation*}
The first part of   condition [\textbf{C}$_3(r)$] means that the vertex $r\bv{g}$ is an MPP for the set~$rG$. Note that under the second part of the condition, this MPP $r\bv{g}$ for  $rG$ is unique (e.g., by Lemma~1 in~\cite{PaBo}). Since  $ \bv{N}(r)$ always belongs to the closure of $Q^+$, the third part of   condition [\textbf{C}$_3(r)$] just excludes the case when the normal $ \bv{N}(r)$ to the level surface of~$\Lambda$ at the point~$r\bv{g}$ belongs to the boundary of the set $rG$.

\begin{theorem*}
Let conditions {\rm [\textbf{C}$_1$], [\textbf{C}$_2$]} and {\rm [\textbf{C}$_3(r_G)$]} be met. Then the asymptotic relation~\eqref{first} holds true,
where $D(G)$ is the value of the second rate function~\eqref{Dv} on $G$ and the constant $A_0 \in (0, \infty)$ can be computed  explicitly.
\end{theorem*}

The value of the constant  $A_0 \in (0, \infty)$ is given by the limit as $\delta\to 0$ of the expressions given by formula~(68) in~\cite{PaBo} for the distribution of~$\bv{\xi}\deq \bv{X}(\delta)$. When proving the theorem below, we demonstrate that that limit does exist and is finite and positive.

\begin{proof}
For a $\delta >0,$ consider the embedded random walk $\{\bv{X}(n\delta)\}_{n\in \N}$ and, for a set $V\subset\R^d$,  denote the first time that random walk hits that set~$V$ by
\[
  \eta_\delta (V) :=\inf\{n\in \N: \bv{X}(n\delta)\in V\}.
\]
First observe that, on the one hand, for any $\delta>0,$ one clearly has
\begin{equation}
\label{>s}
\pr (\tau (sG)<\infty)\ge \pr (\eta_\delta  (sG)<\infty) .
\end{equation}
On the other hand, assuming without loss of generality that $\min_{1\le j\le d}g_j\ge 1$ and setting $I(s):= (\tau(sG), \tau(sG)+\delta]\subset \R$ on the event $\{\tau (sG)<\infty\},$ we have, for any $\ep>0,$
\begin{align}
\pr (\eta_\delta  ((s-\ep)G)<\infty)
% & \ge
% \pr \bigl(\tau (sG)<\infty,\bv{X} (\lceil \eta_\delta (sG)\rceil\delta)\in sG\bigr)
%
 & \ge
   \pr \Bigl(\tau (sG)<\infty,
    \sup_{  t \in I(s)} \|\bv{X} (t) - \bv{X}(\tau (sG))\|\le \ep\Bigr)
   \notag
   \\
 & = \pr  (\tau (sG)<\infty)
    \pr \Bigl(\sup_{ t \in I(s)}\|\bv{X} (t) - \bv{X}(\tau (sG))\|\le \ep\Bigr)
    \notag
    \\
 & = \pr  (\tau (sG)<\infty)
    \pr \Bigl(\sup_{ t \in (0,\delta]}\|\bv{X} (t) \|\le \ep\Bigr),
    \label{rf3}
\end{align}
where the last two relations follow from the strong Markov property and homogeneity of~$\bv{X}$.

Now take  an arbitrary small $\ep>0$. As the process $\bv{X}$ is right-continuous, there exists a $\delta (\ep)>0$ such that
\[
\pr \Bigl(\sup_{ t \in (0,\delta(\ep)]}\|\bv{X} (t) \|\le \ep\Bigr)>(1+\ep)^{-1},
\]
which, together with~\eqref{rf3},  yields  for all $\delta\in (0,\delta(\ep)]$   the inequality
\begin{equation}
\label{<s}
\pr  (\tau (sG)<\infty)\le (1+\ep) \pr (\eta_\delta  ((s-\ep)G)<\infty).
\end{equation}

The precise asymptotics of the probability on the RHS of~\eqref{>s} were obtained in~\cite{PaBo}. It is given in terms of the second rate function $D^{[\delta]}$ for the distribution of the jumps $\bv{X}(n\delta)-\bv{X}((n-1)\delta)\deq\bv{X}(\delta)$ in the random walk $\{\bv{X}(n\delta)\}_{n\ge 0}$. Recalling the well-known fact that the cumulant of $\bv{X}(\delta)$ is given by $\delta K$, we see that the first rate function $ \Lambda^{[\delta]}$ for $\bv{X}(\delta)$ equals
\begin{align*}
%\label{rf1}
  \Lambda^{[\delta]} (\bv{\alpha})  & =
  \sup_{\sbv{\lambda} \in \Theta_{\psi}} (\langle \bv{\alpha},\bv{\lambda} \rangle - \delta K(\bv{\lambda}))
   \\
   &= \delta \sup_{\sbv{\lambda} \in \Theta_{\psi}} (\langle \bv{\alpha}/\delta,\bv{\lambda} \rangle -  K(\bv{\lambda}))
   =\delta \Lambda (\bv{\alpha}/\delta)
   , \quad  \bv{\alpha} \in \R^d
\end{align*}
(cf.~\eqref{rf1}). Therefore the second rate function (see~\eqref{Dv})  $D^{[\delta]}$ for $\bv{X}(\delta)$ is
\[
D^{[\delta]} (\bv{v}) := \inf_{u > 0}  u\Lambda^{[\delta]}( \bv{v}/u)
 = \inf_{u > 0} (u \delta)  \Lambda (\bv{\alpha}/(u \delta)) =D (\bv{v}),
   \quad  \bv{v} \in \R^d.
\]
That is, the second rate function for the random walk  $\{\bv{X}(n\delta)\}$ is the same for all $\delta>0$, which makes perfect sense as one would expect the same asymptotics for the probabilities $\pr (\eta_\delta  (sG)<\infty)$ for different~$\delta$. Hence  the respective value $r^{[\delta]}_G$ (see~\eqref{mpt})  can easily be seen to be given by $\delta r_G$.

%Continuing in a similar fashion to use the above representations for $ \Lambda^{[\delta]}$ and $D^{[\delta]}$ to compute for the random walk $\{\bv{X}(n\delta)\}$ all the values participating in the expression for the constant $A=A^{[\delta]}$ given by~(68) in~\cite{PaBo}, one can verify that it actually does not depend on~$\delta$ as well.

Therefore, applying Theorem~1 in~\cite{PaBo} to the random walk $\{\bv{X}(n\delta)\}$ and using notation~$A^{[\delta]}$ for the constant $A$ appearing in that theorem for the said random walk, we conclude that, for any $\delta\in (0,\delta(\ep)],$  as $s\to \infty,$
\begin{equation}
% \label{lastassertion}
\pr\big( \eta_\delta(sG) < \infty \big) = A^{[\delta]}s^{-(d-1)/2}e^{-sD(G)}(1+o(1)),
\end{equation}
and likewise
\begin{equation}
% \label{lastassertion}
\pr\big( \eta_\delta((s-\ep)G) < \infty \big) = A^{[\delta]} (s-\ep)^{-(d-1)/2}e^{-(s-\ep) D(G)}(1+o(1)) .
\end{equation}
Now from~\eqref{>s} and~\eqref{<s} we see that, as $ s\to
\infty,$
\[
A^{[\delta]}(1+o(1))
\le
R(s):=\frac{\pr (\tau (sG)<\infty)}{ s^{-(d-1)/2}e^{-sD(G)}}
\le
 \frac{A^{[\delta]}(1+\ep)e^{\ep D(G)}}{(1-\ep/s)^{(d-1)/2}}\, (1+o(1))
\]
Therefore, setting $\underline{R}:= \liminf_{s\to\infty}R(s),$ $\overline{R}:= \limsup_{s\to\infty}R(s),$ we have
\[
A^{[\delta]} \le \underline{R}
 \le  \overline{R} \le (1+\ep)e^{\ep D(G)} A^{[\delta]}
\]
for any $\delta\in (0,\delta(\ep)],$ and hence
\[
\limsup_{\delta \to 0} A^{[\delta]}
 \le \underline{R}
 \le  \overline{R}
 \le (1+\ep)e^{\ep D(G)} \liminf_{\delta \to 0} A^{[\delta]}.
\]
As $\ep>0$ is arbitrary small, we conclude that there exists $\lim_{\delta\to 0} A^{[\delta]}=:A_0\in (0,\infty).$ Therefore there also exists
\[
\lim_{s\to\infty}R(s)=A_0.
\]
The theorem is proved.
\end{proof}

{\bf Acknowledgements.} The authors are grateful to MATRIX Research Institute for hosting and supporting the Mathematics of Risk program during which they obtained the result presented in this note. This work was partially supported by Polish National Science Centre Grant  No.~2015/17/B/ST1/01102
(2016-2019) and the ARC  Discovery grant DP150102758.

The authors are also grateful to Enkelejd Hashorva who pointed at a bug in the original version of the note.

\end{document}